\documentclass[final,onefignum,onetabnum]{siuro250211}

\usepackage{lipsum}
\usepackage{booktabs}

\usepackage{amsfonts}
\usepackage{amssymb}
\usepackage{graphicx}
\usepackage{epstopdf}
\usepackage{algorithm}
\usepackage{algpseudocode}
\algrenewcommand\algorithmicrequire{\textbf{Input:}}
\usepackage{float}
\usepackage{xspace}
\usepackage[caption=false]{subfig}

\newcommand{\N}{\mathbb{N}}
\newcommand{\R}{\mathbb{R}}

\renewcommand{\vec}[1]{\boldsymbol{#1}}
\newcommand{\mat}[1]{\mathbf{#1}}

\renewcommand{\d}{\partial}

\newcommand{\pc}{Predictor--Corrector\xspace}
\newcommand{\reopt}{reoptimization\xspace}
\newcommand{\Reopt}{Reoptimization\xspace}
\newcommand{\tb}{\bar{\vec\theta}}
\newcommand{\ttild}{\tilde{\vec\theta}}
\newcommand{\dt}{\Delta\vec{\theta}}

\newcommand{\z}{\vec z}
\newcommand{\zb}{\bar{\vec z}^*}

\newcommand{\zstar}{\vec z^*}

  \ifpdf
  \DeclareGraphicsExtensions{.eps,.pdf,.png,.jpg}
  \else
  \DeclareGraphicsExtensions{.eps}
  \fi

\usepackage{enumitem}
\setlist[enumerate]{leftmargin=.5in}
\setlist[itemize]{leftmargin=.5in}

\title{Efficient Path Continuation for PDE-Constrained Inverse Problems
Under Multiple Perturbations of the Auxiliary Parameters}

\ifdefined\dedication
  \author{
    Carly Haggard\footnote{All authors contributed equally.}
    \thanks{Westmont College (\email{chaggard@westmont.edu}).}
    \and John Introne\footnotemark[2]
    \thanks{William \& Mary (\email{jwintrone@wm.edu}).}
    \and Zoltan Lauko\footnotemark[2]
    \thanks{Macalester College (\email{zlauko@macalester.edu}).}
    \and Jessica Sapozhnikov\footnotemark[2]
    \thanks{Princeton University (\email{js7081@princeton.edu}).}
  }
  
  \dedication{\small\textit{Project advisors: Alen Alexanderian, Joseph Hart, and Troy Larsen}}
\else
  \author{Carly Haggard \and John Introne \and Zoltan Lauko \and Jessica Sapozhnikov}
\fi

\usepackage{amsopn}

\begin{document}

\maketitle

\begin{abstract}
In an inverse problem, we use a mathematical model and state data to reconstruct quantities in the model that are not directly observable. We consider physical systems governed by partial differential equations (PDEs) with unknown quantities to reconstruct and uncertain {\em auxiliary parameters} that remain fixed during inversion. To measure uncertainty in the solution of the inverse problem due to changes in the auxiliary parameters, we solve the inverse problem under multiple perturbations of those parameters. One approach is to use a standard optimization algorithm, such as Newton's method, to solve the inverse problem independently for each perturbation of the auxiliary parameters. However, this reoptimization strategy neglects previous computations and becomes expensive when there are many perturbations, due to the increasing number of required PDE solves. In this work, we develop an efficient method for solving inverse problems under multiple perturbations of the auxiliary parameters. We employ a pseudo-time continuation framework that uses sensitivity information to continuously track the solution across the auxiliary parameter space. This is combined with a {\em minimum spanning tree} (MST) over a set of perturbations. The performance of the proposed method is analyzed under different computational configurations and against similar methods. Numerical experiments indicate that the proposed method reduces the number of required PDE solves by a factor of three compared to a standard reoptimization approach. 
\end{abstract}
\ifdefined\dedication\else
  \noindent\textbf{Author affiliations and contact information}\par
  \noindent Carly Haggard, Westmont College (chaggard@westmont.edu).\par
  \noindent John Introne, William \& Mary (jwintrone@wm.edu).\par
  \noindent Zoltan Lauko, Macalester College (zlauko@macalester.edu).\par
  \noindent Jessica Sapozhnikov, Princeton University (js7081@princeton.edu).\par
  \noindent All authors contributed equally.\par
  \noindent\textit{Project advisors: Alen Alexanderian, Joseph Hart, and Troy Larsen.}\par
\fi

\section{Introduction}
\label{sec:intro}

Inverse problems emerge in science and engineering fields with applications including geophysics, climate science, tomography, and chemical engineering \cite{Antil_Kouri_Lacasse_Ridzal_2018,applications,survey-inverse}. Typically, physical models govern these problems and contain hidden structures that influence state observations, which are the only quantities that we can directly measure. Inverse problems seek to infer these unknown \emph{inversion parameters} from state data. For example, we may want to reconstruct the diffusivity of a material based on data of how heat moves through the material. The system also usually depends on external conditions that we cannot precisely determine. We focus on parameterized inverse problems governed by PDEs, which can be posed as optimization problems of the form
\begin{equation}\label{eq:general}
    \min_{\z} J(\z;\vec{\theta}) 
\end{equation}
for an objective function $J$ that implicitly depends on the PDE. Here, $\z$ consists of the inversion parameters, and $\vec{\theta}$ is a vector of auxiliary parameters that remain fixed during optimization. In applications, $\vec{\theta}$ contains model quantities that are subject to uncertainty which we cannot estimate precisely, such as approximations of physical constants and measurements that incur noise and systematic error.

Ideally, the inverse problem would jointly optimize the inversion and auxiliary parameters, but the problem becomes prohibitively expensive in high-dimensional settings with many auxiliary parameters. Computing gradients and Hessians of the objective requires multiple PDE solves, which are computationally intensive. By fixing the auxiliary parameters, solving the inverse problem becomes manageable. It also requires us to quantify the uncertainty in the solution due to changes in the fixed auxiliary parameters. That is, we must repeatedly solve the inverse problem under perturbations of the auxiliary parameters to measure uncertainty in the inversion parameters, which may be expensive depending on the optimization algorithm.

We can further reduce the computational cost of solving the inverse problem in two ways: by simplifying the problem constraints or improving the algorithm's efficiency. First, model reduction improves the tractability of solving inverse problems by finding low-rank approximations of the PDE map and other operators, as well as constraining the state and parameter spaces \cite{martinez-reduction,BOYCE2014168,lieberman-reduction}.
Model reduction increases uncertainty in the solution resulting from changes to the auxiliary parameters. Between fixing the auxiliary parameters and simplifying the model, computing the sensitivity of the solution to changes in the parameters becomes essential to understanding the behavior of the true inversion quantities. In this paper, we explore the second way of reducing computational cost: improving the optimization algorithm's efficiency. Specifically, we examine how to decrease the number of derivative (especially Hessian) calculations in each iteration of the optimization algorithm and the total number of required iterations to converge to the solutions under multiple parameter perturbations.

Given a minimizer of the objective for a nominal auxiliary parameter vector, we aim to efficiently compute the minimizers for perturbed parameter vectors. The natural idea to solve inverse problems across perturbations is \emph{\reopt}, which entails running a standard optimization algorithm such as Newton's method for each perturbation. In \reopt, the initial iterate of each run can be a previously computed solution, but it largely treats each optimization problem as independent of one another. This method is inefficient in that it ignores information about how the minimizer changes across the parameter space. As the inverse problem becomes more challenging and high-dimensional, \reopt may require an excessive number of total iterations to converge and PDE solves at each iteration. We assert that considering the sensitivity of the solution to the auxiliary parameters, reusing computations, and optimally traversing the parameter space reduces computational cost.

In this work, we develop a more efficient method that employs a pseudo-time continuation framework \cite{hart2026preconditioned}. The idea is to continuously track the solution to the inverse problem along paths in the parameter space. At each point in the path, the algorithm begins by computing the sensitivity of the minimizer to the auxiliary parameters to estimate the solution at a perturbed parameter vector through a linear approximation. Unless the minimizer function is linear within the neighborhood of the step, this first order prediction incurs some error. To account for this error, the method performs corrections using a standard optimization algorithm to recover the true minimizer. There are two primary advantages of this framework over \reopt. Given a reasonably sized step along the parameter path, the linear approximation closely estimates the true minimizer curve without the high cost of an iteration from a standard optimization algorithm. Secondly, the method reuses derivative information to speed up calculations of the objective's Hessian. 

This algorithm solves the inverse problem for a single parameter perturbation. To find the solutions for multiple perturbations, we must determine an optimal path to traverse the parameter space. This context leads to the combination of tools from graph theory with pseudo-time continuation. We treat the auxiliary parameter vectors as vertices of a graph, where directed edges represent solving the inverse problem for terminal vertices starting from the solutions of the outgoing vertices. One option is repeated \reopt, where we run the optimization algorithm independently for each vertex, except the initial iterates are the solutions from nearby vertices. Prior work on continuation found \reopt to be inefficient for single perturbations \cite{hart2026preconditioned}, but it remains unclear how \reopt performs in the setting of multiple perturbations with varying distances from the nominal point.

Our contributions in this paper are threefold. We fuse pseudo-time continuation with \emph{minimum spanning trees} (MSTs) to solve inverse problems under multiple auxiliary parameter perturbations, highlighting improvements in efficiency over \reopt. We construct the MSTs using Prim's algorithm \cite{Prim} along with an edge weight function to estimate the expense of performing continuation along each edge. Secondly, we implement a Hessian preconditioner with persisting state throughout continuations runs, which further enhances efficiency. Finally, we compare three distinct edge weights in constructing the MST. The rest of the paper is structured as follows. Section~\ref{sec:prelim} outlines pseudo-time continuation. Section~\ref{sec:method} describes our algorithm that extends continuation for multiple perturbations. Section~\ref{sec:formulation} establishes the model inverse problem that we analyze. Section~\ref{sec:results} provides results of our method for the model inverse problem. We conclude our discussion and include lines of future work in Section~\ref{sec:conclusion}.

\section{Preliminaries}
\label{sec:prelim}

We aim to efficiently solve parameterized PDE-constrained inverse problems. Defining the weighted norm as $\|\z\|^2_\mat{R} = \z^\intercal \mat{R} \z$ for a symmetric positive definite matrix $\mat{R}$, we pose these inverse problems as optimization tasks of the form
\begin{equation}
\label{eq:objective}
\min_{\z} \; J(\z; \vec{\theta}) := 
  \tfrac{1}{2} \| \mathcal{F}(\z; \vec{\theta}) - \vec{y} \|^2 + \tfrac{1}{2} \|\z \|_\mat{R}^2,
\end{equation} 
where $\mathcal{F}$ represents the discretized forward PDE map that depends on the inversion parameters $\z$ and auxiliary parameters $\vec{\theta}$, $\vec{y}$ consists of the state measurement data, and $\mat{R}$ is a regularization matrix. The solution $\zstar$ to  Equation~\eqref{eq:objective} depends on $\boldsymbol{\theta}$, which remains fixed during optimization, and each minimizer satisfies the first-order optimality condition
\begin{equation}\label{eq:optimality}
    \nabla_{\z} J(\zstar;\vec{\theta})=\vec{0}.
\end{equation}
We assume that $J$ is twice continuously differentiable with respect to $\z$ and $\vec{\theta}$. Also, we assume that the Hessian $\mat{H} = \nabla_{\z\z} J$ is invertible in the region of the parameter space where we wish to solve the optimization problem as we must be able to apply its inverse to compute the minimizer's sensitivity to the parameters. In most practical settings, both assumptions hold. We may then apply the implicit function theorem to Equation~\eqref{eq:optimality} to write the solution $\zstar(\vec{\theta})$ as a differentiable function of the auxiliary parameters.

\subsection{\pc}
\label{sec:pc}

To efficiently solve the inverse problem under a perturbation of the auxiliary parameters, we employ a pseudo-time continuation framework that consists of steps along linear approximations to the minimizer function and corrections to these approximations that recover the true minimizers. We refer to this approach as \pc. Suppose that we have found the minimizer $\zb$ for the auxiliary parameter vector $\tb$ and aim to solve the inverse problem for a single perturbation $\ttild$. Letting $\dt = \ttild - \tb$, our goal is to continuously track the minimizer as the parameters vary along the linear path 
\begin{equation}\label{eq:path}
    \vec{\theta}(t) = \tb + t\dt, \qquad t \in [0,1].
\end{equation} 
Let $\mat{D} = \frac{\d\zstar}{\d\boldsymbol{\theta}}$ denote the Jacobian of $\zstar(\vec{\theta})$. Differentiating the first-order optimality condition \eqref{eq:optimality} by $\vec{\theta}$  and evaluating each component at $(\zstar(\vec{\theta}), \vec{\theta})$ gives the relation $\mat{H} \mat D = - \mat{B}$, where $\mat{B} = \nabla_{\z\vec{\theta}}J$ is the mixed second derivative of the objective. Along $\vec{\theta}(t)$, let the minimizer curve be $\zstar(t) = \zstar(\vec{\theta}(t))$, which exists locally and is differentiable when $\mat{H}$ is invertible \cite{fiacco1983sensitivity}. We have the following differential equation for the minimizer:
\begin{equation}\label{eq:ode}
    \frac{d\zstar}{dt} = \mat{D}(\zstar(t),\vec\theta(t))\dt, \qquad \zstar(0) = \zb.
\end{equation}
We traverse $\vec{\theta}(t)$ in $N$ equal-sized iterations of size $\Delta t = 1/N$. For $k \in \{0,\dots,N\}$, we denote $\vec\theta_k = \vec\theta(k\Delta t)$ and $\zstar_k$ the minimizer at $\vec\theta_k$. At each iteration, \pc performs a linear approximation (the predictor step) followed by iterations of a standard optimization algorithm (the corrector steps) if necessary. In this work, the predictor utilizes the first-order forward Euler method for the differential equation in \eqref{eq:ode}: 
\begin{equation}\label{eq:predictor}
    \zstar_{\mathrm{pred}} = \zstar_k + \Delta t \mat{D}_k\dt,
\end{equation}
where the subscript $k$ denotes evaluation at $(\zstar_k,\vec{\theta}_k)$. Higher-order Runge--Kutta methods are also viable for the predictor. The greater accuracy of these methods may allow for larger step sizes, but there is a tradeoff with extra sensitivity computations. The true minimizer function is often nonlinear, so this linear predictor step introduces error that we must correct before taking another predictor step. In this work, the corrector executes the Inexact Newton-Conjugate Gradient (CG) \cite{nocedal2006numerical} method with $\zstar_{\mathrm{pred}}$ as the initial iterate:
\begin{equation}\label{eq:corrector}
    \mat{H}(\zstar_{\mathrm{pred}}, \vec{\theta}_{k+1})\,\vec{p} = \nabla_{\z}J(\zstar_{\mathrm{pred}};\vec{\theta}_{k+1}),
    \qquad
    \zstar_{k+1} = \zstar_{\mathrm{pred}} - \eta\,\vec{p},
\end{equation}
where $\eta$ is chosen by a line search with Wolfe conditions \cite{wolfe}. CG is an iterative method for solving linear systems involving a symmetric positive definite (SPD) matrix \cite{iterative-linear-systems}. Directly solving or inverting these systems is expensive in high dimensions, but CG often reaches an accurate approximation within few iterations if the matrix is well-conditioned or has tight clusters of eigenvalues. However, note as an exception that inexact Newton-CG does not require the Hessian to be positive definite. If the algorithm encounters a direction of negative curvature in CG, the previous iterate is still a direction in which the objective decreases \cite{nocedal2006numerical}. The method is inexact since $\vec{p}$ may not perfectly solve the system in \eqref{eq:corrector} at every iteration.

\subsection{Crucial Components of \pc}
\label{sec:components}

Since the differential equation in \eqref{eq:ode} tracks the minimizer, the predictor steps require that each $\zstar_k$ is a solution to the inverse problem. To account for error in the predictor at each iteration, we repeat corrector steps until the gradient of the objective reaches a small fixed tolerance. Assuming that the step size is not too large, correcting until tolerance ensures that \pc returns solutions of the same accuracy (that is, with roughly the same gradient norm) as \reopt. Thus, differences in computational costs between methods arise solely from differences in efficiency.

The Hessian provides one of the largest challenges to efficiency. In high-dimensional ill-posed inverse problems, the Hessian is often ill-conditioned. Ill-conditioning causes numerical errors and increases the required number of CG iterations to solve linear systems to some tolerance. Since \pc involves inverting the Hessian at each step, the ill-conditioning increases computational costs. To combat this issue, constructing a reliable preconditioner accelerates CG convergence \cite{iterative-linear-systems}. For \pc, the preconditioner uses information from the objective's curvature and past CG computations to reduce the condition number of the Hessian \cite{hart2026preconditioned}. Additionally, building and storing matrices is intractable, so in implementations, we only have access to the action of matrices like the Hessian on vectors.

When closely tracking the minimizer, the Hessian remains SPD in most practical settings, but it may become indefinite if the step size is too large or the predictor is inaccurate. While Inexact Newton-CG can handle an indefinite Hessian, the method may require many iterations to converge. Thus, choosing the right step size is crucial to the success of \pc. Unfortunately, the optimal $N$ varies in each problem due to differences in the nonlinearity of the minimizer function. A smaller step size admits better linear approximations of the minimizer trajectory and therefore performs fewer corrector steps but requires more iterations to traverse the path $\vec{\theta}(t)$. Conversely, large steps tend to make the predictor inaccurate, requiring additional corrector steps. The step size generally only affects the cost of traversing the parameter space since the corrector continues until a gradient tolerance. We assume that the step size is small enough so that the predictor does not push the estimated solution outside of the initial Newton basin. Otherwise, \pc may converge to a different local minimizer than \reopt since solutions are not unique in ill-posed inverse problems.

\section{Method}\label{sec:method}

\pc provides a framework for tracking the solution to an inverse problem under a single perturbation of the auxiliary parameters. This section describes our approach that extends \pc to the context of multiple perturbations. Section~\ref{sec:sampling} outlines our approach to generating a set of sample perturbations. Section~\ref{sec:mst_construction} details the construction of the MST for the set of perturbations. Section~\ref{sec:algorithm}  provides pseudocode for a simplified version of the algorithm, omitting the details of the preconditioner. Section~\ref{sec:comp_consider} concludes by highlighting additional algorithmic considerations used to increase efficiency.

\subsection{Auxiliary Parameter Sampling}
\label{sec:sampling}

Suppose that we have solved the inverse problem for a nominal point $\tb$, which in practice is an estimate of the unknown ground truth auxiliary parameters. To account for uncertainty in this estimate, we solve the inverse problem for many auxiliary parameter perturbations sampled around $\tb$. We define a probability distribution on the auxiliary parameter space and sample a set $\Theta = \{\ttild_1, \dots, \ttild_\nu\}$ of parameter vectors from that probability distribution. In this work, we use a Gaussian distribution centered at $\tb$ with covariance elements related to fractional margins of error in the auxiliary parameters.

\subsection{MST Construction} 
\label{sec:mst_construction}

We must build a path in the parameter space that contains $\tb$ and every point in $\Theta$ before employing \pc to solve the inverse problem for each vector. We represent the perturbations as a weighted digraph $G = (V,E)$, where the vertices $V$ are auxiliary parameter vectors, $E \subseteq V \times V$ represents directed edges between vertices, and each edge has a weight determined by a function $w: V\times V\to \R^{\ge 0}$ (not necessarily a metric). We must form edges between the vertices to create a spanning tree. Here, edges refer to running \pc from the initial to final vertices. In this work, $V = \Theta \cup \{\tb\}$ and $G$ is complete, meaning $E$ consists of all ordered pairs between distinct vertices. As a baseline, we consider what we deem the \emph{naive approach}, which simply selects the edges that connect each perturbation to the nominal point $\tb$. This method creates a spanning tree with a large total path length in the parameter space. A more efficient idea constructs an MST with respect to the weight $w(\vec{\theta}_1, \vec{\theta}_2)$ defined for each edge.

To build the MST, we use Prim's greedy algorithm \cite{Prim}. Starting at $\tb$, we add the edge of lowest weight to the MST and append the connecting vertex to the list of visited vertices. Iteratively, we choose the edges of lowest weight among those that are connected to the visited vertices until we have a spanning tree. The algorithm naturally returns an order of edges outward from $\tb$ which guarantees that we have computed the minimizers at starting vertices when performing continuation. Generally, the cost to build the MST is much smaller than repeatedly solving the inverse problem when $\nu$ is not too large, so having $\mathcal{O}(\nu^2)$ edges is not a concern. Reducing the cost to construct the MST may become more important as the number of vertices significantly increases, which we discuss in the conclusion. Traversing from $\vec{\theta}_1$ to $\vec{\theta}_2$, the simplest weight is Euclidean distance:
\begin{equation}
\label{eq:euclid}
    w_E(\vec{\theta}_1, \vec{\theta}_2) = \| \vec{\theta}_1 - \vec{\theta}_2 \|^2 = \| \Delta \vec{\theta} \|^2.
\end{equation}
Since this weight is symmetric, $G$ simplifies to a complete undirected graph.

\begin{figure}[H]
    \centering
    \includegraphics[width=0.9\linewidth]{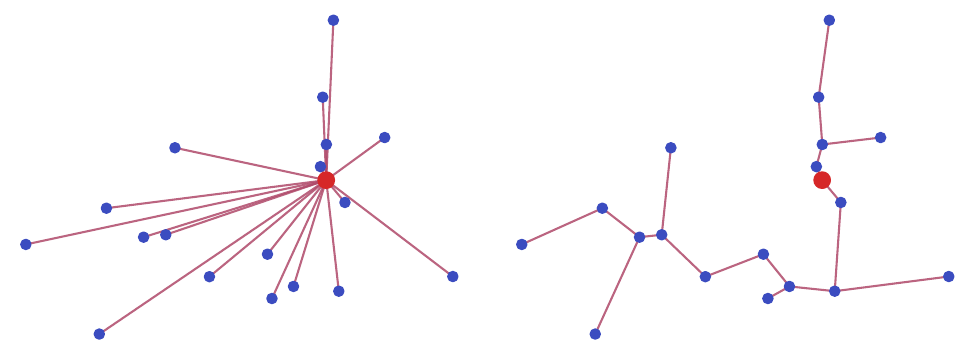}
    \caption{Left: The spanning tree using the naive approach, which branches from $\tb$ (red vertex) to each $\ttild_i$ (blue vertices). Right: The Euclidean MST over the same set of perturbations.}
    \label{fig:clouds}
\end{figure}

Figure~\ref{fig:clouds} depicts the different spanning trees using the naive approach and an MST with Euclidean distance for the same $\Theta$. Many edges in the naive tree are longer than those in the MST tree, and the total path length is much shorter for the MST tree. The one benefit of the naive approach is that it does not require selection of an edge weight, but our results in Section~\ref{sec:results} indicate that using $w_E$ is consistently better.

\subsection{Algorithm}
\label{sec:algorithm}

We present the key details of our method in Algorithm~\ref{alg:method}. The primary inputs are the number of perturbations $\nu$ to sample from the parameter distribution, the edge weight function $w$ in constructing the MST, and the number of steps $N$ to take along each edge of the MST. After generating the set of perturbations $\Theta$ and building the MST with Prim's algorithm, we employ \pc along every edge of the tree in order of the returned list of edges from Prim's algorithm.

Each run of \pc follows the outline from Section~\ref{sec:prelim}. We compute the first derivative of the minimizer with respect to the auxiliary parameters by solving the linear system involving the Hessian and mixed second derivative of $J$. The predictor approximates the new solution by extending the previous minimizer in the direction of the derivative according to the step size. Then, we perform Newton corrector iterations until the objective gradient norm reaches a small fixed tolerance. If the predictor is very accurate, no corrector steps are necessary. Typically, using a good step size will mean that only one or two corrections are necessary. If the step size is too large, the method will likely require many corrections. We solve linear systems involving the Hessian with CG. To accelerate CG convergence, we use a preconditioner, which we discuss in detail in Section~\ref{sec:prec}.

\begin{algorithm}
\caption{\pc for multiple perturbations.}
\label{alg:method}
\begin{algorithmic}[1]
    \Require Number of perturbations $\nu \in \N$, edge weight function $w$, number of steps $N \in \N$
    \State $\Delta t \gets 1/N$
    % \State Initialize $\mat{P} = \mat{R}^{-1}$ for $\tb$
    \State Generate $\nu$ parameter samples $\Theta = \{\ttild_1, \dots, \ttild_\nu\}$
    \State Build MST using Prim's algorithm with $w$
    \For{edge in MST from $\vec{\theta}_1$ to $\vec{\theta}_2$}
        \State $\Delta \vec{\theta} \gets \vec{\theta}_2 - \vec{\theta}_1$
        \State Initialize $\zstar_0$ with the minimizer for $\vec{\theta}_1$
        \For{$k = 0,\dots, N-1$}
            \State Solve for $\mat{D}_k$ in $\mat{H}_k \mat{D}_k = -\mat{B}_k$ using CG \label{alg:method:cg1}
            \State $\zstar_{\text{pred}} \gets \zstar_k + \Delta t \mat{D}_k \Delta\vec{\theta}$
            \While{$\|\nabla_{\z} J\| \ge \text{tolerance}$}
                \State Solve for $\vec{p}$ in $\mat{H} \vec{p} = \nabla_{\z} J$ using CG \label{alg:method:cg2}
                \State Update $\zstar_{k+1}$ in the direction of $-\vec{p}$
            \EndWhile
        \EndFor
    \EndFor
    \State \Return Final minimizers $\zstar$ for each $\ttild_i$
\end{algorithmic}
\end{algorithm}

\subsection{Computational Considerations}
\label{sec:comp_consider}

There are several opportunities to further improve the computational efficiency and robustness of our method in the context of multiple perturbations by optimizing how we traverse the parameter space.

\subsubsection{Alternative Edge Weights}
\label{sec:costs}

While Euclidean distance is simple and effective, it becomes less informative in high dimensions and is agnostic to structure in the inverse problem. A better edge weight considers how the minimizer evolves across parameter space.

Two options that we analyze are sensitivity informed and physics informed weights:
\begin{equation}
\label{eq:alt_costs}
\begin{aligned}
    w_S(\vec{\theta}_1, \vec{\theta}_2)
    &= \| \mat{D}(\zstar_1, \vec{\theta}_1) \Delta\vec{\theta} \|^2, \\
    w_P(\vec{\theta}_1, \vec{\theta}_2)
    &= \left\|
        \mathcal{F}(\zstar_1; \vec{\theta}_1)
        - \mathcal{F}(\zstar_1; \vec{\theta}_2)
       \right\|^2,
\end{aligned}
\end{equation}
where $\mat{D} = \frac{\d \zstar}{\d \vec{\theta}}$ is the Jacobian of the minimizer function $\zstar(\vec{\theta})$ introduced in Section~\ref{sec:pc}, $\zstar_1$ is the minimizer corresponding to $\vec{\theta}_1$, and
$\mathcal{F}$ is the PDE map in Equation~\eqref{eq:objective}. The intuition behind the sensitivity weight is that if the minimizer does not change rapidly in the direction of a perturbation, the predictor error may be small relative to what Newton's method can quickly correct. Then fewer corrector iterations are necessary, and a larger step size may be more feasible. The minimizer function may still be highly nonlinear, so this weight may overestimate the benefit of traveling in directions with small first derivative. A more informative weight would incorporate higher derivatives, but deriving those tensors through a closed form expression or finite difference calculations would incur high overhead costs or approximation errors. The idea behind the physics weight is that if the PDE solutions under two different auxiliary parameter vectors are similar when evaluated at the same diffusion $\zstar_1$, then the solution $\zstar_2$ at $\vec{\theta}_2$ is likely similar to $\zstar_1$. We then expect \pc to be more efficient because the predictor is more accurate.

Finally, in practice, we must make simplifications to these weights for them to be practical. As noted in Section~\ref{sec:components}, building the matrix $\mat{D}$ is expensive. Attempting to compute $\mat{D}$ at every parameter vector requires more PDE solves than actually running \pc over $\Theta$. Instead, we build the full matrix at the nominal point $\tb$ and reuse it for all perturbations. Similarly, we replace $\zstar_1$ with $\zb$ in the physics weight, so we only solve the PDE once for the nominal point and each perturbation to compute all edge weights. That is, MST construction demands $\nu+1$ PDE solves instead of $\mathcal{O}(\nu^2)$. Consequently, we can form a more informed MST without significant overhead costs. Note that Euclidean distance has a major advantage in this regard as there is no overhead cost to derive the MST. Since $\Theta$ consists of points within a fairly small neighborhood of $\tb$, extrapolating information from the nominal point leads to small losses in accuracy. It also allows us to simplify $G$ to be undirected and construct the MST upfront instead of iteratively as we solve the inverse problem for more perturbations.

\subsubsection{Preconditioner}
\label{sec:prec}

Preconditioning the Hessian significantly enhances convergence of CG in lines~\ref{alg:method:cg1}~and~\ref{alg:method:cg2} of Algorithm~\ref{alg:method}. We employ the adaptive preconditioner from previous continuation work \cite{hart2026preconditioned}. We initialize it as the inverse of the regularization matrix. After running CG during the predictor step in line~\ref{alg:method:cg1}, we perform a block quasi-Newton update of the preconditioner, similar to BFGS methods \cite{bfgs}, using previously computed CG iterates. After the first iteration of the corrector, we again improve the preconditioner with a block quasi-Newton update and also perform a parametric quasi-Newton update that approximates the objective's curvature from gradient calculations. The computational cost of \pc lies mostly in solving the Hessian for each CG iteration. With this adaptive preconditioner, CG becomes cheaper after just one step along the edge. Thus, for high-dimensional inverse problems, taking small steps avoids the risk of the Hessian becoming indefinite without substantially increasing the computational cost over using the optimal step size.

Updating the adaptive preconditioner naturally extends to the context of multiple perturbations. We initialize the preconditioner at the base point $\tb$ as the inverse of the regularization matrix as before. After traversing each edge, the final state of the preconditioner contains information that is still relevant for a subsequent \pc run that starts at that terminal vertex. As such, when traversing the set of perturbations $\Theta$, vertices inherit the state of the preconditioner instead of resetting it to the inverse of the regularization matrix.

\subsubsection{Step Size Adaptivity}
\label{sec:adapt}

A key component to efficiency, robustness, and practicality of this algorithm lies in the step size. For a small model problem, evenly-spaced steps are effective, and we can investigate performance across values of the number of steps $N$. In a highly nonlinear problem, the optimal step size will change throughout the path $\vec{\theta}(t)$, and the algorithm must handle unexpected conditions. While adaptivity is the focus of future work, we pose two ideas. First, as discussed in Section~\ref{sec:components}, \pc is inefficient when the Hessian becomes indefinite. Preliminary analyses indicate that by backtracking and reducing the step size when the Hessian becomes indefinite in CG, the algorithm automatically corrects an initial step size that is too large with minimal overhead. This idea provides a mechanism to shrink the step size. A complementary idea is to increase the step size if it seems that we could travel a larger distance in parameter space without increasing the number of required CG and Newton iterations (e.g., when the predictor error is near zero) \cite{hackl1980efficient}. Together, the system should efficiently and robustly traverse the parameter space without prior knowledge of the optimal step size, allowing for practical use. In this work, we initialize $N=2$, which we found to provide consistently fast performance.

\section{Example PDE}\label{sec:formulation}

We analyze an infinite-dimensional inverse problem for a steady-state advection-diffusion PDE similar to the inverse problem from  \cite{alexanderian2023perspective}: 
\begin{equation}\label{eq:strong}
    \begin{cases}
        -\big(\kappa(x)u'(x)\big)' + \theta(x)u'(x) = s(x), & x \in (0,1), \\
        \kappa(x)u'(x) = \alpha u(x), & x = 0, \\
        \kappa(x)u'(x) = -\alpha u(x), & x = 1.
    \end{cases}
\end{equation}
Advection transports the entire state $u(x)$ with a spatially varying wind velocity field $\theta(x)$, which we take to be the auxiliary parameter. Diffusion evenly disperses the state throughout the domain. The spatially varying diffusion $\kappa(x)$ is the unknown quantity to reconstruct. The inversion parameter $z(x)$ is the log-diffusion (i.e., $\kappa(x) = e^{z(x)}$), which ensures that the estimated diffusion is positive. The source term $s(x) = \exp\{-50(x-\frac{1}{2})^2\}$ is fixed and independent of $\kappa$ and $\theta$. In the Robin boundary conditions, the coefficient is $\alpha = 1$. 

\subsection{Inverse Problem Formulation}
\label{sec:inverse_prob}

Since the diffusion to estimate is infinite-dimensional, we first pose the inverse problem with continuous components. The objective is
\begin{equation}\label{eq:opt}
\min_{z(x)} \; J(z(x);\theta(x)) = 
  \frac{1}{2}\int_0^1 \big( u(x) - y(x) \big)^2 \, dx + \frac{\beta}{2}\,\mathcal{P}(z),
\end{equation}
where $y(x)$ is the observed state data, $\beta > 0$ is the regularization
coefficient, and $\mathcal{P}(z)$ is the regularization operator for the inversion parameter defined as
\begin{equation}\label{eq:reg}
    \mathcal{P}(z) = \int_0^1 z(x)^{2} \, dx + 10^{-3} \int_0^1 z'(x)^{2} \, dx.
\end{equation}
The true log-diffusion and advection field in our experiments are given by
\begin{equation}\label{eq:true}
    z_\text{true}(x) = 1 + 2\sin^2 (2\pi x), \qquad \theta_{\text{true}}(x) = 1 + 0.3(1-x)^2.
\end{equation}
Together, they generate the true state, and we measure  state observations with $2\%$ multiplicative Gaussian noise. The inverse problem is ill-posed, which usually means that the solution is non-unique. Solving the PDE uniquely determines $u(x)$ given $z(x)$ and $\theta(x)$, so we pose the problem in reduced space with $z(x)$ as the only optimization variable. Evaluating $J$ requires one PDE solve, and evaluating its gradient requires an additional adjoint solve \cite{adjoints}. Each Hessian-vector product requires two further solves. We measure computational cost as the total number of PDE solves since they are the most expensive operation in the optimization.

\subsection{Regularization}\label{sec:reg}

The first term in Equation~\eqref{eq:opt} is a data misfit that measures the discrepancy between the predicted and observed state. While the prediction is ideally very close to the observed data, without regularization, the reconstruction is sensitive to small changes and noise in the observed data. The regularization in Equation~\eqref{eq:reg} includes $L_2$ penalties on the magnitude and spatial variation of the log-diffusion $z(x)$, which encourage stable and smooth reconstruction of the solution for this ill-posed problem. We empirically determined the $10^{-3}$ coefficient on the spatial variation term to be best for reconstructing the diffusion, and $\beta$ provides an overall weight for regularization in the objective.

The value of $\beta$ dramatically affects both the accuracy and efficiency of solutions to the inverse problem. If too small, the minimizer overfits to noise in the data and poorly matches the true parameter. The objective and minimizer function become less well-behaved, so both \reopt and especially \pc become slower. Conversely, when $\beta$ is too large, the solution becomes over-constrained and cannot meaningfully track variation in the true inversion parameter. The functions become smooth enough that the predictor accrues minimal error. In either case, the reconstructed inversion parameter is not accurate.

One way to determine $\beta$ invokes Morozov's discrepancy principle, which holds that the data misfit should have the same magnitude as the noise in the state data \cite{NAIR_2009}. Not only does this method require prior knowledge of the noise level, but we also found that the resulting regularization is too strong. Another approach uses the prior in a Bayesian formulation of the inverse problem \cite{bui-thanh}. In this work, we utilize the L-curve method, which locates the point of highest curvature on a log-log plot of solution norm as a function of the data misfit \cite{lcurve}. The regularization at this point balances the complexity of the solution with its fit to the data. We validated that the L-curve's ideal $\beta$ of $2.33 \times 10^{-7}$ that we use closely matches the $\beta$ which minimizes the difference between the optimal diffusion and ground truth in this problem.

\subsection{Discretization}\label{sec:discretize}

We discretize the PDE in \eqref{eq:strong} using continuous piecewise-linear finite elements on a uniform mesh of $n = 100$ nodes with spacing $h = 1/(n-1)$. This process yields a mass matrix $\mat{M}$ and stiffness matrix $\mat{S}$. The mass-weighted inner product $\z^{\intercal}\mat{M}\z$ represents the squared $L_{2}$ norm of the log-diffusion on $(0,1)$ and the stiffness-weighted inner product $\z^{\intercal}\mat{S}\z$ the squared $L_{2}$ norm of its derivative. We weight the norms from previous sections with these matrices. We denote $\z$ and $\vec{\theta}$ for the discretized vectors of $z(x)$ and $\theta(x)$, $\vec{u}$ for the corresponding PDE solution, and $\vec{y}$ for the vector of measurements of the true state. Writing $\mat{R} = \mat{M} + 10^{-3} \mat{S}$, the discretized objective in the form of Equation~\eqref{eq:objective} becomes
\begin{equation}\label{eq:discretized_obj}
    \min_{\z} \; J(\z;\vec{\theta}) = 
  \tfrac{1}{2} \| \vec{u} - \vec{y} \|_\mat{M}^2 + \tfrac{\beta}{2}  \|\z \|^2_\mat{R}.
\end{equation}

\begin{figure}
    \centering
    \includegraphics[width=1\linewidth]{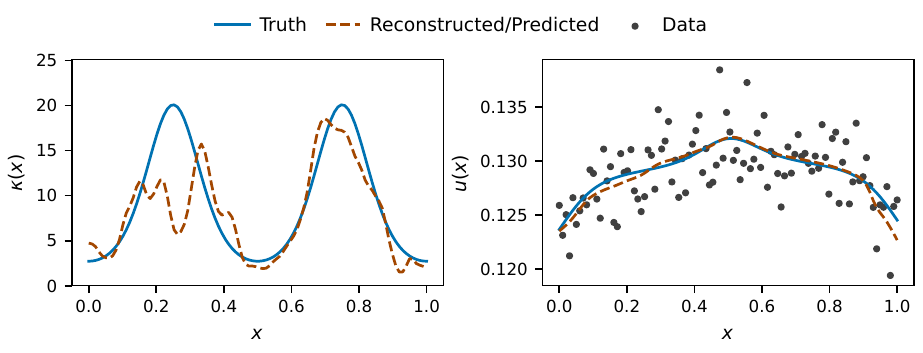}
    \caption{Left: Reconstructed diffusion coefficient $\kappa(x)$. Right: Corresponding predicted state $u(x)$. The inverse problem is solved at the nominal auxiliary parameter $\tb$.}
    \label{fig:inverse_fit}
\end{figure}

The nominal auxiliary parameter vector $\tb$ comes from the function $\bar{\theta}(x) = 1 + (1-x)^2$. Notice that the coefficient is $1$ instead of $0.3$ for $\theta_{\text{true}}$ in \eqref{eq:true}, so $\tb$ overestimates the strength of the true advection field. Before generating a set of perturbations, we use Newton's method to solve the inverse problem at this base point. Figure~\ref{fig:inverse_fit} depicts the differences between the ground truth diffusion and state and their estimates under $\tb$. While the optimized diffusion generally matches the shape of the truth, it fails to capture the peak for $x$ around $0.2$ to $0.3$. This error results from solving the inverse problem under $\tb$ instead of the ground truth and using noisy observation data in an ill-posed problem.

\subsubsection{Perturbations}\label{sec:perturbations}

\begin{figure}
    \centering
    \includegraphics[width=1\linewidth]{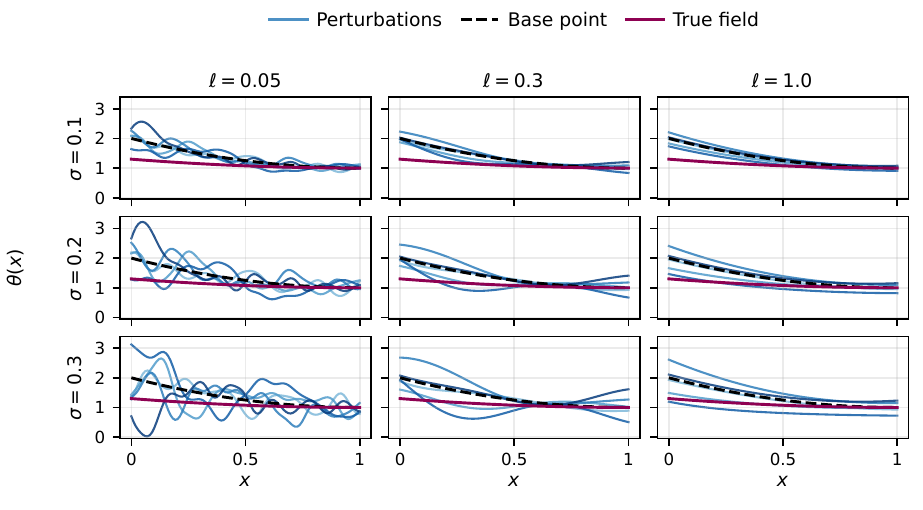}
    \caption{Sample perturbations about $\tb$. Each panel shows 5 perturbed advection fields (blue), the nominal field $\tb$ (black) and the true advection field $\boldsymbol{\theta_\text{true}}$ (magenta) for different values of $\sigma$ and $\ell$.}
    \label{fig:perturbatios}
\end{figure}

To sample parameter perturbations about $\tb$, we construct the distribution $\mathcal{N}(\vec{1}, \mat{\Sigma})$ and multiply the entries element-wise to $\tb$, where
\begin{equation}
\label{eq:kernel}
\Sigma_{ij} = \sigma^2 \exp\left( - \frac{ \| x_i - x_j \|^2}{2\ell^2} \right).
\end{equation}
Here, $\sigma$ represents the fractional margin of error in the entries of $\tb$ and $\ell$ the length scale. As $\ell$ approaches $0$, the sample adds independent white noise to each coordinate. At larger values, the noise in each coordinate transitions more smoothly. We identify $\ell$ as a measure of correlation among the coordinates of each perturbation. Figure~\ref{fig:perturbatios} depicts the effect of different $\sigma$ and $\ell$ on sample perturbations. As $\sigma$ increases, the amplitude of noise about $\tb$ increases. As $\ell$ increases, samples fluctuate less and coordinates are either all greater or all less than those of $\tb$. By solving the inverse problem across these samples, we can quantify the uncertainty in the minimizer over a space that (roughly) contains the ground truth. In our experiments, we test up to $\sigma=0.3$. The range of $\ell$ from $0.05$ to $1.0$ captures a sufficiently wide range of smoothness in the parameters to understand performance of the method.

\section{Computational Results}\label{sec:results}

We conduct three primary studies on the efficiency of our method for the inverse problem defined in Section~\ref{sec:formulation}. First, we measure the benefit of inheriting the preconditioner. Then we compare our MST approach to \reopt and naive continuation. Finally, we compare the performance for different MST edge weight functions. In all experiments, we measure computational cost as the total number of PDE solves required to solve the inverse problem for all perturbations.

\subsection{Preconditioner Inheritance}
\label{sec:prec_results}
\begin{figure}
    \centering
    \includegraphics[width=0.9\linewidth]{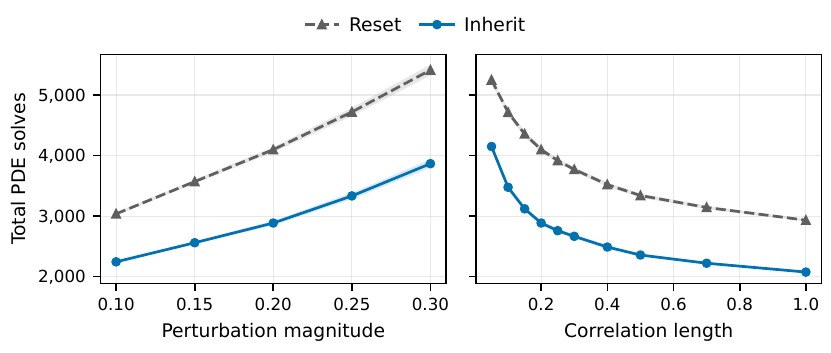}
    \caption{Computational cost of inheriting versus resetting the preconditioner throughout solving each perturbation with $95\%$ confidence intervals. Left: Across $\sigma$ with $\ell=0.2$. Right: Across $\ell$ with $\sigma=0.2$.}
    \label{fig:inheritance}
\end{figure}

To verify our discussion in Section~\ref{sec:prec}, we determine whether inheriting the preconditioner throughout the tree of perturbations decreases overall computational cost over resetting the preconditioner for each continuation run. In these experiments, we perform $50$ independent runs using the Euclidean MST for $\nu=100$ perturbations and $N=2$ steps. We empirically found two steps to be consistently efficient for this problem. Figure~\ref{fig:inheritance} depicts the results when varying the perturbation magnitude parameter $\sigma$ and correlation parameter $\ell$. In both cases, inheriting the preconditioner provides a meaningful reduction in the number of PDE solves, with an average $28\%$ improvement across both perturbation magnitude and correlation. We also found that the improvement is consistent across the depth of the vertices in the tree. For vertices that are three or more edges away from the base point, inheriting the preconditioner cuts about $25\%$--$35\%$ of the total cost. With these results, we use inheritance for the MST methods in the rest of our experiments.

\subsection{\Reopt vs. Continuation}
\label{sec:reopt_vs_cont}

Previous work found that \pc provides significant performance gains over \reopt for single perturbations \cite{hart2026preconditioned}. We test if this trend extends to the regime of multiple perturbations. We compare the computational cost of two baselines, \reopt and the naive approach, to our Euclidean MST. To present the strongest results for \reopt, we use a ``warm-started" version where the initial iterate for each perturbation is the solution of the closest parameter vector \cite{gondzio1998warm}. Essentially, this method takes advantage of the Euclidean MST, except it uses Newton's method instead of \pc. That is, it does not consider minimizer sensitivity information, nor does it use an adaptive preconditioner. As the distance between the starting and target perturbations increases, the initial iterate becomes less accurate, and more Newton iterations are necessary to solve the inverse problem. The naive approach uses \pc, but every edge begins at $\tb$. Using a fixed number of continuation steps, distant vertices cause large predictor errors and intermediate estimates that may land outside of the Newton basin. Reusing solutions in the proposed MST approach provides for shorter edges, which mitigates the predictor error and reduces the number of corrector iterations. Building the Euclidean MST requires no extra PDE solves, and vertices inherit the preconditioner.

\begin{figure}
    \centering
    \includegraphics[width=1\linewidth]{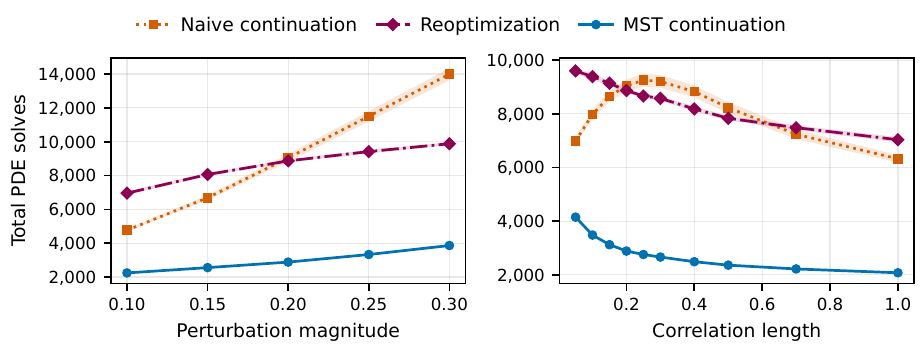}
    \caption{Computational cost of \reopt versus the naive approach versus the Euclidean MST with $95\%$ confidence intervals. Left: Across $\sigma$ with $\ell = 0.2$. Right: Across $\ell$ with $\sigma = 0.2$.}
    \label{fig:reopt_vs_cont}
\end{figure}

We randomly sample $\nu=100$ perturbations for $50$ independent runs and use $N=2$ steps for the continuation methods. Figure~\ref{fig:reopt_vs_cont} displays the efficiency of the three methods across perturbation magnitudes and correlation levels. Notably, the Euclidean MST requires on average approximately one-third as many PDE solves as \reopt across perturbation magnitude and length scale. As $\sigma$ increases, all methods require more PDE solves, reaching $13{,}991$ for the naive approach and $9{,}882$ for \reopt while only $3{,}864$ for the Euclidean MST at $\sigma = 0.30$. The costs for the continuation approaches rise steadily, whereas the cost for \reopt slightly levels off. The Euclidean MST is also much faster across $\ell$. At $\ell = 1.0$, it requires about $2{,}000$ solves compared to over $6{,}000$ for the naive approach and $7{,}000$ for \reopt. The Euclidean MST continuation and \reopt see a consistent decrease in cost as $\ell$ increases, but the naive approach surprisingly becomes slower for moderate levels of correlation before surpassing \reopt again at larger $\ell$.

The two baselines are more expensive for different reasons. \Reopt has no predictor or adaptive preconditioner. The cost levels off at larger $\sigma$ because as the initial iterate becomes less accurate, not many more iterations are necessary to converge. The cost decreases across $\ell$ because its warm-started iterates take advantage of the decreasing average distance between perturbations. Naive continuation does have a predictor and preconditioning, but since the edges all begin at $\tb$, the predictor error grows quadratically in the edge length for fixed $N$. Naive continuation is therefore the cheaper baseline for small perturbations but more expensive beyond $\sigma = 0.20$. Unlike the other two methods, the naive approach does not consistently improve as $\ell$ increases. For this particular PDE, the relative change in the solution compared to that of the base point increases for moderate levels of correlation with the same shape as the cost displayed in the figure. Continuation with the MST retains the advantage of the predictor without concerns of strong decreases in accuracy for large perturbations or moderate levels of correlation by keeping parameter paths short.

\subsection{MST Edge Weight Functions}
\label{sec:cost_results}

The results from Section~\ref{sec:reopt_vs_cont} indicate that pairing an MST with \pc is the most efficient approach, but is the Euclidean weight the best choice for this problem? We conduct similar experiments to Section~\ref{sec:reopt_vs_cont} but compare the cost when using an MST with the Euclidean weight against MSTs with the two weights from Equations~\eqref{eq:alt_costs}: sensitivity informed, and physics informed. As discussed in Section~\ref{sec:costs}, we simplify these two alternative weight functions by considering information from the base point to make the MST construction tractable. All three weights use the same parameter samples, regularization, and continuation settings for $50$ trials of $\nu=100$ perturbations. The reported PDE solve counts include both MST construction and continuation. 

\begin{figure}
    \centering
    \includegraphics[width=1\linewidth]{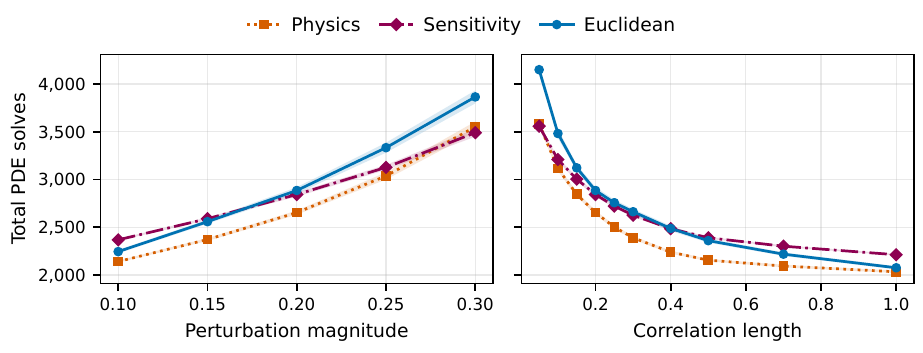}
    \caption{Computational cost of MST continuation using Euclidean, sensitivity informed, and physics informed weights with $95\%$ confidence intervals. Left: Across $\sigma$ with $\ell = 0.2$. Right: Across $\ell$ with $\sigma = 0.2$.}
    \label{fig:alternative_costs}
\end{figure}

Figure~\ref{fig:alternative_costs} depicts the computational cost of the three weights across perturbation magnitude and correlation. All three have the same general trends, but the two alternative weights highlight improvements over Euclidean distance. The physics informed weight requires fewer PDE solves than the Euclidean weight at every tested perturbation magnitude and correlation level, reducing the total PDE solves by $7.6\%$ on average in the perturbation magnitude experiments and $9.2\%$ on average in the correlation level experiments. Using the sensitivity informed weight is more expensive than the other two for small perturbations and large correlations levels. Nonetheless, it is better than Euclidean distance in many of the configurations and the best of the three for $\sigma = 0.3$ in the left plot and $\ell = 0.05$ in the right plot.

In fact, the sensitivity weight reduces the number of PDE solves during continuation the most but has a steep overhead in MST construction. Building the full sensitivity matrix $\mat{D}$ at the base point required about $402$ PDE solves per trial, whereas the physics informed weight required $101$ solves. This additional overhead outweighs the savings during continuation for small perturbation magnitudes and large correlation lengths, where the Euclidean weight is cheaper, but has less of an impact for the more difficult experimental settings. Also, the cost to build $\mat{D}$ increases steadily with the problem dimension but does not depend on the number of perturbations. Conversely, the overhead for the physics informed weight scales with the number of perturbations but does not change with the dimension. Thus, the performance of an MST with a better informed edge weight depends on whether the overhead to construct the MST is small relative to the reduction in continuation costs, and that occurs in different settings for the different weights. Note that using a fixed step size for \pc, we did not realize the full benefit of the sensitivity weight in choosing the best step size.

\section{Conclusion}
\label{sec:conclusion}

We have presented an efficient method for solving parameterized PDE-constrained inverse problems under multiple perturbations of the auxiliary parameters. The method constructs an MST with respect to a user-defined edge weight function over the set of samples. Starting from the nominal point, the tree determines the optimal piecewise-linear path to follow and solve the inverse problem for each perturbation. \pc then tracks the minimizer along each edge by iteratively computing the sensitivity of the solution to the auxiliary parameters and correcting until reaching a gradient tolerance.

For a model advection-diffusion PDE, we determined an efficient configuration for this method. We found that inheriting the adaptive preconditioner between nodes decreases computational cost over simply resetting it for every continuation run. Across perturbation magnitude and length scale, the method outperformed the two baselines, requiring about one-third as many PDE solves as warm-started \reopt. The predictor step and adaptive preconditioner allow continuation to run faster than \reopt, and the MST mitigates the distance to travel in parameter space, which protects the predictor's accuracy under a fixed number of steps. While Euclidean distance is a useful metric, other MST edge weights such as the physics or sensitivity informed weights provide further efficiency gains.

%\subsection{Limitations and Future Work}

We tested our method on a small steady state PDE with one spatial dimension and a single spatially-varying auxiliary parameter. We intend to test the method on a two-dimensional problem that is more computationally challenging due to a higher dimensionality and more complicated parameter mesh. Previous work suggests that continuation becomes more powerful over \reopt as the problem becomes more challenging \cite{hart2026preconditioned}, so we will analyze how that trend extends for multiple perturbations.

Another limitation is that we generated data on the same mesh that we used for inversion, and every mesh coordinate provided an observation. Ideally, the two meshes would be different, and data would come from an observation operator or real measurements. We held the base point $\tb$ fixed without varying its distance from $\vec{\theta}_{\text{true}}$, which may affect how efficient continuation is over \reopt. We also held the number of continuation steps constant at $N = 2$. Determining the best step size is a key challenge for this method. While preconditioning allows us to take conservatively sized steps without much extra cost, the step size is pivotal for the efficiency and robustness of the method. It could also change which MST edge weight is best. Future work will investigate whether an adaptive step size like the ones in \cite{hackl1980efficient,KUBATSCHEK2025107747} or Runge--Kutta methods \cite{DORMAND198019} lead to stronger performance alongside backtracking.

Additionally, we only sampled perturbations using the outlined Gaussian distribution. The ultimate goal of this work is to efficiently and accurately quantify the uncertainty and sensitivity of the solution to the inverse problem under changes to the auxiliary parameters. To appropriately model the auxiliary parameters in practice, the underlying distribution may be more complex and require changes to the spanning tree construction. First, if the distances between perturbations is large, adding intermediate vertices may reduce the total cost to traverse the MST. One idea to this end involves clustering the perturbations and adding the centroids to the graph. If the number of perturbations significantly increases, the computational cost to construct the MST may be a concern when edges connect every pair of vertices. We can then use the clusters to reduce the number of edges by only connecting all vertices within each cluster and extending a single connection between clusters.

Secondly, constructing a more sophisticated, nonlinear path in the parameter space may enhance efficiency. For example, if the perturbations tightly follow a nonlinear curve, we may fit a smoothing spline \cite{green1993nonparametric} through the points and create offshoots to each perturbation. We may also design a nonlinear path around how the minimizer evolves across parameter space. Similar to the sensitivity weight, we could iteratively travel in a direction of small minimizer sensitivity. While a piecewise-linear path is likely suitable in most cases, these ideas may be useful in specific problems. Notably, we saw that different MST weights were optimal under different sampling settings even for the same problem. We assert that the three proposed weights should be studied further in addition to ideas like a statistically informed weight such as Mahalanobis distance that captures information about the parameter distribution.

Finally, extensions of this method include a higher-order \pc approach \cite{higher_pc} and asymptotic numerical methods \cite{COCHELIN19941181}. Instead of taking linear approximations to the minimizer function under changes to auxiliary parameters, these methods track nonlinear changes. Generally, fewer continuation steps are necessary. However, computing higher-order sensitivities also increases the cost of each step. Since the predictor estimate simply needs to be relatively close to the true solution for the corrector to quickly converge, a first-order predictor is sufficient in most problems. Nonetheless, higher-order methods are worth investigating as they may improve efficiency in highly nonlinear problems.

\section*{Acknowledgments}
This work was supported by the National Science Foundation under awards DMS\#2349611 and DMS\#2342344. The authors thank their mentors for their continued support. 

The authors used Claude Code, Google Gemini, and ChatGPT to assist with ideation and code implementation. The authors assume responsibility for all content. 

\bibliographystyle{siamplain}
\bibliography{references}

\end{document}